\documentclass[a4paper,12pt]{article}
\usepackage{graphicx}
\usepackage{amssymb}
\usepackage{amsmath, amscd, amsthm}
\usepackage{mathrsfs}
\usepackage[all]{xy}
\usepackage{amsfonts}
\usepackage{amsrefs}
\usepackage{tikz}
\usepackage{tikz-cd}
\begin{document}

\theoremstyle{remark}
\newtheorem*{Remark}{Remark}
\newtheorem*{Remarks}{Remarks}
\newtheorem*{Proof}{Proof}
\theoremstyle{definition}
\newtheorem{Definition}{Definition}
\newtheorem*{Example}{Example}
\newtheorem*{Examples}{Examples}

\theoremstyle{plain}
\newtheorem{Theorem}{Theorem}
\newtheorem{Theorema}{Theorem}
\newtheorem{Proposition}{Proposition}
\newtheorem{Corollary}{Corollary}
\newtheorem{Lemma}{Lemma}

\newcommand{\QQ}{\mathbb Q}

\renewcommand{\labelenumi}{\roman{enumi})}
\newcommand{\vs}{\vspace{5mm}}
\renewcommand{\baselinestretch}{1.10}
\setlength{\parskip}{3mm}
\setlength{\parindent}{0em}
\small\normalsize

\newcommand{\C}{k}
\newcommand{\MM}{\mathfrak M}
\newcommand{\M}{\mathbb M}

\begin{center}
{\Large\sc A motivation for some local cohomologies 
}\\ $\ $ \\
Ricardo Garc\'{\i}a L\'opez
\footnote{
 Departament de Matem\`atiques i Inform\`atica.
 Universitat de Barcelona (UB),
 Gran Via, 585.
 E-08007, Barcelona, Spain.
 e-mail: ricardogarcia@ub.edu
}
\footnote{\today}
\end{center}

These notes are an expanded version of my talk at the FACARD workshop
in Barcelona, January 2019. Comments on them are  wellcomed. 
\vs

{\bf 0. Introduction }

If $k$ is a field, $R$ is a commutative Noetherian $k$-algebra and $I\subset R$ is an ideal, the local cohomology modules 
$H^i_I(R)$ have been fruitfully studied by Gennady Lyubeznik through the consideration of supplementary structures on them, besides the structure of $R$-module: A 
$\mathcal D$-module structure if $\text{char}(k)=0$, a $F$-module structure if $\text{char}(p)=p>0$ (see \cite{Lyu}, \cite{Lyu2}). These structures are functorial, in the sense that
if $I\subset J$ are ideals and $i\geqslant 0$, they are preserved by the natural maps $H^i_J(R) \longrightarrow H^i_I(R)$, and they are also preserved by connecting homomorphisms in long exact sequences and by the differentials in the Mayer-Vietoris spectral sequence for local cohomology. 

One might ask the following, somewhat imprecise question: What is ``the richest additional structure'' one can put on such local cohomologies so that functoriality is preserved? It turns out that, in a certain sense, some of the theorems proved by Madhav Nori and Florian Ivorra for the study of (perverse) motives of algebraic varieties can be applied in this context, and they provide a rather precise, albeit quite abstract, answer. The purpose of these notes is to expose briefly the theorems of Nori and Ivorra and to show in which way they are related to the 
question about local cohomology sketched above. We also discuss on ``motivic'' versions of the Lyubeznik numbers, a
set of numerical invariants for local rings introduced in \cite{Lyu}. Two short digressions are included at the end.
\medskip

{\bf 1. Nori motives }

Motives were envisioned by Grothendieck as an attempt to find a common source for cohomology theories of algebraic varieties. 
If $\mathfrak{Pman}$ is the category of projective manifolds, say over $\mathbb C$, we have several functors
\[
H_{sing}, H_{DR}, H_{\ell-adic}, H_{cris},.. : \mathfrak{Pman} \longrightarrow \QQ, \mathbb R, \QQ_p,..-\text{vector spaces}\ ,
\]
and the conjectural category of pure motives should be an abelian, $\mathbb Q$-linear, monoidal category $\mathfrak{Mot}$\footnote{Usually one expects additional properties: Tannakian, with finite dimensional Homs,...} endowed with
a contravariant ``universal'' cohomology functor $\mathfrak{Pman} \longrightarrow \mathfrak{Mot}$ so that the above functors (and more generally, any Weil cohomology theory) factor through it. For varieties which are possibly non-proper or singular, the conjectural category of mixed motives should have analogous properties\footnote{See for instance the introduction to \cite{Lev}*{Part I} for a list of expected properties.}.
At present, the existence of such categories seems far from being settled (although already Grothendieck proposed a conjectural category of pure motives,
Voevodsky constructed a category with properties one would expect for the derived category of motives, and there are other proposals, see e.g. \cite{Lev}, \cite{And} and the references therein). 
\medskip

Madhav Nori introduced in some lectures\footnote{To my knowledge, he himself never published anything on this subject.} a new approach which leads to the construction of a candidate for the category of mixed motives. For the rest of this section, we follow \cite{HMS}:

\underline{Idea:} Use the singular cohomology functor to 
construct a category which is universal {\it for cohomology theories functorially comparable with singular cohomology}.\footnote{More generally, one can fix a cohomology theory 
to define ``a candidate for the category of motives'', and then this category is universal for cohomology theories functorially comparable to the chosen one.} It goes roughly as follows:

A diagram is an oriented graph (possibly infinite). 
A category $\mathcal C$ gives rise to a diagram that we still denote by $\mathcal C$, or sometimes by $D( \mathcal C)$ (vertices are the objects of the category, edges are the morphisms, we forget about composition).
If $D$ is a diagram and $\mathcal C$ a category, a representation of $D$ is a morphism of oriented graphs $D \longrightarrow \mathcal C$.
If $k$ is a field, we denote by $k-\mathfrak{fvs}$ the category of finite dimensional $k$-vector spaces. 

\begin{Theorem} (Nori's theorem)
Let $D$ be a diagram, $T: D \longrightarrow k-\mathfrak{fvs}$ a representation. There is a $k$-linear abelian category $\mathcal C(D,T)$,
a representation
\[
\widetilde{T}: D \longrightarrow \mathcal C(D,T) \ ,
\]
and a faithful, exact, $k$-linear functor $\mathcal F_T: \mathcal C(D,T) \longrightarrow k-\mathfrak{fvs}$ such that $T= \mathcal F_T \circ \widetilde{T}$.

Moreover, if $\mathcal A$ is another $k$-linear abelian category endowed with a representation $U: D \longrightarrow \mathcal A$ and a faithful, $k$-linear exact functor $\mathcal F_U: \mathcal A \longrightarrow k-\mathfrak{fvs}$ with $T= \mathcal F_U \circ U$, then there is a faithful, $k$-linear functor $\mathcal C(D,T) \longrightarrow \mathcal A$ making the following diagram commutative\footnote{Commutativity of diagrams of functors should be understood up to canonical natural transformations, see \cite{Ivo} for a more precise statement.}
\[
\xymatrix{ {} & 
{\mathcal C(D,T)} \ar@/_0.1 pc/[1,1]^-{\mathcal F_T}  \ar@/^0.1pc/[2,0]
& {} \\ {D } \ar@/^0pc/[-1,1]^-{ \widetilde{T}}
\ar@/^0.1pc/[0,2]_-{\ \ \ \ \ T}
\ar@/_0.1pc/[1,1]_-{U}
& {} 
& { k-\mathfrak{fvs}.} \\ {} 
& {\mathcal A } 
\ar@/_0pc/[-1,1]_-{\mathcal F_U}   & {}}
\]
The category $\mathcal C(D,T)$ together with $\widetilde{T}, \mathcal F_D$ are uniquely determined by this property up to unique equivalence of categories. 
 $\mathcal C(D,T)$ will be called a diagram category or a Nori category.
\end{Theorem}

\underline{About the proof:} Nori defines the $k$-algebra $\text{End}_k(T)$ of endomorphisms of $T$ as the subalgebra of 
\[
\prod_{v \text{ vertex of D}} \text{End}_k(Tv)
\]
formed by those $(\varphi_v)_v$, $\varphi_v\in\text{End}_k(Tv)$, such that for every edge $e$ from $v$ to $w$ the square
\[ \begin{CD}
Tv @>Te>> Tw\\
@V\varphi_vVV @VV\varphi_wV\\
Tv @>Te>> Tw
\end{CD} \]
is commutative. If $D$ is a finite diagram, then $\mathcal C(D,T)$ is the category $\text{End}_k(T)-\mathfrak{mod}$ of those left $\text{End}_k(T)$-modules which are finite dimensional as $k$-vector spaces.
In the general case, consider the set $\mathfrak F(D)$ of finite full subdiagrams of $D$. Then $\mathcal C(D,T)$ is defined as a colimit (more precisely, a $2$-colimit) of the categories $\text{End}_k(T_{\mid F})-\mathfrak{mod}$ where $F$ runs over $\mathfrak F(D)$.  

In more detail: Given $F_1, F_2\in \mathfrak F(D)$ if $F_1\subset F_2$, there is a natural morphism of $k$-algebras
$\text{End}_k(T_{\mid F_2})\longrightarrow \text{End}_k(T_{\mid F_1})$ (given an element of $\text{End}_k(T_{\mid F_2})$, consider only its components corresponding to vertices in $F_1$), and so a restriction of scalars functor $\mathcal R_{F_1,F_2}: \mathcal C(F_1,T_{\mid F_1})\longrightarrow\mathcal C(F_2,T_{\mid F_2})$. The objects of $\mathcal C(D,T)$ are those of the categories $\mathcal C(F,T_{\mid F})$ where $F\in \mathfrak F(D)$, and where an object $X$ in $\mathcal C(F_1,T_{\mid F_1})$ is identified with $\mathcal R_{F_1,F_2}(X)$ for any $F_2\supset F_1$. If $X,Y$ are objects, then $X$ is an object in some $\mathcal C(F_1,T_{\mid F_1})$ and $Y$ in some $\mathcal C(F_2,T_{\mid F_2})$, and the morphisms from $X$ to $Y$ are
\[
 Hom_{\mathcal C(D,T)}(X,Y) = \varinjlim_{F_3\supset F_1,F_2} Hom_{\mathcal C(F_3,T_{\mid F_3})}(\mathcal R_{F_1,F_3}(X),\mathcal R_{F_2,F_3}(Y))
\]
The functor $\mathcal F_T:\mathcal C(D,T) \longrightarrow k-\mathfrak{fvs}$ is the forgetful functor. $\spadesuit$
\medskip

\begin{Remarks}
 \begin{enumerate}
  \item In general, $\mathcal C(D,T)$ is not equivalent to the category of modules over an algebra. But it is equivalent to the category $A(D,T)-\mathfrak{comod}$ of comodules over the coalgebra $A(D,T):=\varinjlim_{F}\text{End}(T_{\mid F})^{\vee}$ which are finite dimensional over $k$\footnote{The symbol ``$\vee$'' denotes the $k$-dual. One has $A(D,T)^{\vee}=End(T)$, but in general $End(T)^{\vee}\neq A(D,T)$. See digression 1 for more details.}.
  \item Notice that there is complete freedom on choosing the diagram $D$, as long as we have a representation in the category of finite dimensional vector spaces.
 \end{enumerate}
\end{Remarks}
\medskip

\underline{Examples of objects of $\mathcal C(D,T)$:}\ \ Consider for example a vertex $v$ of the diagram $D$. If $F\in \mathfrak F(D)$, then $\text{End}_k(T_{\mid F})$ acts on $Tv$ via the projection $\text{End}_k(T_{\mid F})\longrightarrow \text{End}_k(Tv)$, thus $Tv$ is a $\text{End}_k(T_{\mid F})$-module. 
Then, the construction of $\mathcal C(D,T)$ shows that $Tv$ determines an object of $\mathcal C(D,T)$, which is $\widetilde{T}v$. In fact, every object of $\mathcal C(D,T)$ is a subquotient of a finite direct sum of objects of the form $\widetilde{T}v$ (\cite{HMS}*{Proposition 7.1.16}).
\medskip

The main application of the above construction is obtained when the diagram $D$ is Nori's diagram $\mathfrak{Pairs}$, defined as follows:

\underline{Vertices}: Triples $(X,Y,i)$ where $X$ is a complex quasi-projective variety, $Y\subset X$ a closed subvariety, $i\in\mathbb N$.

\underline{Edges}: There are two types of edges:
\begin{enumerate}
\item \underline{Functoriality edges:} For every morphism $f: X \longrightarrow X'$ with $f(Y)\subset Y'$, there is an edge
\[
(X',Y',i) \longrightarrow (X,Y,i)
\]
\item \underline{Connecting homomorphism edges:} For every triple $X\supset Y\supset Z$, there is an edge
\[
(Y,Z,i) \longrightarrow (X,Y,i+1).
\]
\end{enumerate}

There is a representation:
\begin{eqnarray*}
Sing: \mathfrak{Pairs} &\longrightarrow & \mathbb Q-\mathfrak{fvs}\\
(X,Y,i) &\longmapsto & H_{sing}^i(X,Y; \mathbb Q) \ ,
\end{eqnarray*}
where $H_{sing}^i$ denotes singular cohomology. 
The category of cohomological Nori motives\footnote{To be precise, rational, effective, cohomological Nori motives. One can consider also integer coefficients, homological motives, varieties defined over a subfield of $\mathbb C$, etc...} is the diagram category $\mathcal C(\mathfrak{Pairs}, Sing)$. It has realizations in any cohomology theory comparable with singular cohomology (De Rham, \'etale,...), as an example
we give some details for De Rham cohomology:

\begin{Example}
 Let $Sing_{\mathbb R}: 
\mathfrak{Pairs} \longrightarrow \mathbb R-\mathfrak{fvs}$ denote the representation given by singular cohomology with real coefficients, $\mathcal C(\mathfrak{Pairs}, Sing_{\mathbb R})$ the corresponding Nori category, $DR:\mathfrak{Pairs} \longrightarrow \mathbb R -\mathfrak{fvs}$ the representation given by De Rham cohomology.
As it is well-known, there is an isomorphism of functors between $Sing_{\mathbb R}$ and $DR$. 
So we have
\[
\xymatrix{ {} & 
{\mathcal C(\mathfrak{Pairs},Sing_{\mathbb R})} \ar@/_0.1 pc/[1,1] \ar@/^0.1pc/[2,0]
& {} \\ { \mathfrak{Pairs}} \ar@/^0pc/[-1,1]^-{ \widetilde{Sing_{\mathbb R}}}
\ar@/^0.1pc/[0,2]_-{\!\!\!\!\!\!\!\!\!\!\!\!\!\!\!\!\!\!\!\!\!\!\!\!\!\!\!\!\!\!\!\!\!\!\!\!  Sing_{\mathbb R}}
\ar@/_0.1pc/[1,1]_-{DR}
& {} 
& { \mathbb R-\mathfrak{fvs}.} \\ {} 
& {\mathbb R-\mathfrak{fvs} } 
\ar@/_0pc/[-1,1]_-{\sim}   & {}}
\]
where the vertical arrow exists by the universal property of Nori's category. Hence the De Rham functor factors through $\mathcal C(\mathfrak{Pairs},Sing_{\mathbb R})$\footnote{Also, there is a canonical equivalence $\mathcal C(\mathfrak Pairs, Sing)\otimes \mathbb R \cong \mathcal C(\mathfrak Pairs, Sing_{\mathbb R})$, where for a given $\mathbb Q$-linear abelian category $\mathcal C$,
one denotes $\mathcal C \otimes \mathbb R$ the category which has the same objects as $\mathcal C$ and the morphisms are $Hom_{\mathcal C \otimes \mathbb R}(\cdot,\cdot)=Hom_{\mathcal C}(\cdot,\cdot)\otimes \mathbb R$. In fact,  Nori's construction of a diagram category is compatible with field extensions.}.
\end{Example}
\medskip

Nori's approach is well related to others: 
\begin{enumerate}
\item  Andr\'e introduced in \cite{And2} a definition of (pure) motives, see also \cite{Ara}*{6.4}, \cite{HMS}*{6.1}. Arapura proved that the category of Andr\'e motives is equivalent to the category of pure Nori motives (\cite{Ara}*{Theorem 6.4.1}). A further relation between Andr\'e and Nori motives is proved in \cite{HMS}*{Proposition 10.2.1}. 
 \item There is a category of Nori $n$-motives (it is the thick abelian subcategory generated by the motives corresponding to vertices
 $(X,Y,i)$ with $i\leqslant n$). It has been proved that, for $i=0,1$, they agree up to equivalence with some previously introduced categories (roughly, the category of Nori $n$-motives is equivalent to that of Artin motives for $n=0$, and to the category of Deligne $1$-motives for $n=1$, see \cite{AB}). 
 \item To the category of Nori motives one can attach a ``motivic Galois group'' (one has to introduce a tensor structure, then localize to obtain a neutral Tannakian category). On the other hand, Ayoub introduced also a motivic Galois group for Voevodsky motives. It has been shown by Choudhury and Gallauer that these two groups are isomorphic, see \cite{CG}.
 \item There is a realization functor from the category of Voevodsky effective motives to the derived category of Nori motives (the proof was sketched by Nori, a complete proof is given in the thesis of Harrer, \cite{Harr}). Nori motives have been used for the study of periods, see \cite{HMS}. 
\end{enumerate}

{\bf 2. Perverse motives}

Florian Ivorra generalized Nori's theorem in \cite{Ivo}. 

\underline{\it Idea:} Replace the category of finite dimensional $k$-vector spaces by a essentially small $k$-linear category $\mathcal C$ such that every object is of finite length (following Gabriel, we say that $\mathcal C$ is finite) and with finite dimensional Homs (one says $\mathcal C$ is Hom-finite). The category of perverse sheaves on a complex manifold $X$, or the category of mixed Hodge modules on $X$ are finite and Hom-finite (\cite{BBD}*{Th\'eor\`eme 4.3.1.i)}). Indeed, Ivorra proves that Nori's theorem is also valid if $k-\mathfrak{fvs}$ is replaced by a finite, Hom-finite category. 
\medskip

\begin{Remarks}
\begin{enumerate}
 \item In Ivorra's generalization the obtained diagram category is  equivalent to a category of finite dimensional modules over a coalgebra. An interesting result proved also by Ivorra (\cite{Ivo}*{Corollary 4.3}, partially relying on work by Gabriel and Takeuchi), is that if $k$ is a field, a $k$-linear category $\mathcal C$ is finite and Hom-finite if and only if there is a $k$-coalgebra $A_{\mathcal C}$ and an exact equivalence between $\mathcal C$ and the category $A_{\mathcal C}-\mathfrak{comod}$ of finitely dimensional left $A_{\mathcal C}$-comodules\footnote{A explicit description
of $A_{\mathcal C}$ seems difficult.}. 

\item Nori's theorem has been generalized still much more, see \cite{BCL}*{Remark 2.6(b)}. The methods in loc. cit. allow to extend Nori's theorem to the case where $k-\mathfrak{fvs}$ is replaced by very general categories (including for example abelian categories), they are completely different to those of Nori and Ivorra and rely on mathematical logic. Also, Nori categories can be constructed using Freyd's notion of the abelian hull of a category, see  \cite{BP} or \cite{MI}*{section 1} \footnote{In these generalizations, in general the resulting diagram category is not equivalent to the category of finite dimensional comodules over a $k$-coalgebra.}. 
\end{enumerate}
\end{Remarks}

Given a complex quasi-projective algebraic variety $X$, Ivorra defines the following diagram $\mathfrak D^{Nori}_X$: 

\underline{Vertices}: Triples $(Y \stackrel{a}\to X,Z,i)$ where $Y$ is a complex quasi-projective variety 
and $Z\subset Y$ is closed. Sometimes such a triple will be denoted simply $(Y,Z,i)$.

\underline{Edges}: They are the the following ones:
\begin{enumerate}
 \item Let $(Y_1\stackrel{a_1}\to X,Z_1,i),\, (Y_2\stackrel{a_2}\to X,Z_2,i)$ be triples. If there is a morphism $f:Y_2 \longrightarrow Y_1$ with $a_1\circ f=a_2$ and $f(Z_2)\subset Z_1$, then there is an
 edge 
 \[
  (Y_2\stackrel{a_2}\to X,Z_2,i) \longrightarrow (Y_1\stackrel{a_1}\to X,Z_1,i). 
 \]
\item If $W\subset Z$ is a closed subvariety, then there is an edge\footnote{The indexing is slightly different from the one in \cite{Ivo}.}
\[
 (Y\stackrel{a}\to X,Z,i) \longrightarrow (Z\stackrel{a_{\mid Z}}\to X,W, i+1)
\]
\end{enumerate}

To obtain a representation he considers the category of (rational) perverse sheaves, for which we refer to \cite{BBD}. 
Let $D^b_{cons}(X)$ denote the full subcategory of the bounded derived category of complexes of sheaves of $\mathbb Q$-vector spaces on $X$ whose objects are complexes with constructible cohomology. On $D^b_{cons}(X)$ there is a six functor formalism ($f^{\ast},f^{!}, f_{\ast}, f_{!}$, $\otimes$, $RHom$)
and also a $t$-structure given by the middle perversity, which in turn gives a cohomological functor
\[
^{p}H^0: D^b_{cons}(X) \longrightarrow Perv_{\mathbb Q}(X)
\]
to the heart of this $t$-structure (the abelian category of perverse sheaves, see \cite{BBD}\footnote{An object $K^{\bullet}$ in $D^b_{cons}(X)$ is perverse if 
\[
\dim \text{Supp}\, H^i(K^{\bullet})\leqslant -i \ \ \text{and} \ \ \dim \text{Supp}\, H^{i}(\mathbf{D}K^{\bullet})\leqslant -i \ \ \text{for } \ i\in\mathbb Z
\]
where $\mathbf{D}$ denotes the Verdier dual. Notice that the cohomology of a perverse sheaf is concentrated in non-positive degrees.
}).

Then, he defines a representation $T: \mathfrak D^{Nori}_X \longrightarrow Perv_{\mathbb Q}(X)$ by
\[
 (Y,Z,i) \longmapsto  {^pH}^{i}(a_{!}\,u_{\ast}\,u^{\ast}\,a^{!}(\mathbb Q_X)).
\]
Here, $u: Y\smallsetminus Z \hookrightarrow Y$ is the inclusion, $a^{!},a_{!}, u_{\ast}, u^{\ast}$ are the operations in $D^b_{cons}(X)$ and $^pH$ denotes perverse cohomology. Ivorra defines the category of (effective) perverse Nori motives as the diagram category $\mathcal C(\mathfrak D^{Nori}_X, T)$. 
\medskip

\begin{Remarks}
\begin{enumerate}
\item  When $X=\text{Spec}\,\mathbb C$, one recovers the homological version of Nori's motives in the previous section.

 \item Ivorra considers also Saito's category of mixed Hodge modules (\cite{Sai_MHM}). The category of perverse Nori motives he defines has realizations both in the category of rational perverse sheaves and in the category of mixed Hodge modules (the former factors through the later). 
 \item There is another, perhaps related proposal of D. Arapura of ``motivic constructible sheaves'', based also on Nori's theorem, see \cite{Ara}. Ivorra indicates the possibility of putting a $t$-structure on the derived category of the category of perverse Nori motives so that its heart is equivalent to Arapura's category. Arapura proves that, in his setting, there are direct and inverse image functors, the proof for direct images is difficult. Such a functoriality is not settled in \cite{Ivo} (but see the very recent preprint \cite{MI}).
 \item Ivorra proves that if $X$ is a quasi-projective scheme, then the bounded derived category of perverse Nori motives is related to Ayoub's triangulated category $ D^{et}_{ct} (X,\mathbb Q)$ of \'etale constructible motives with rational coefficients, see \cite{Ivo2}.
%

\item Ivorra considers rational perverse sheaves, we will consider complex perverse sheaves for the rest of this talk. By the Riemann-Hilbert correspondence, if $X$ is a complex algebraic variety\footnote{If $X$ is singular, see digression 2.} of dimension $d_X$, the De Rham functor induces an equivalence of categories 
\begin{eqnarray*}DR:\
 rh\mathcal D_{X}-\mathfrak{mod}  &\longrightarrow &Perv_{\mathbb C}(X) \\
 M\  & \longmapsto & RHom_{\mathcal D_X}(M, \mathcal O_X^{an})[d_X]
\end{eqnarray*}
where $rh\mathcal D_{X}-\mathfrak{mod}$ denotes the category of regular holonomic $\mathcal D_X$-modules and  $M$ is sent to its De Rham complex, and so we can equivalently consider regular holonomic algebraic $\mathcal D_X$-modules.
More precisely, we denote 
$T_{\mathcal D}: \mathfrak D^{Nori}_X \longrightarrow rh\mathcal D_{X}-\mathfrak{mod}$ the representation
\[
T_{\mathcal D}: (Y,Z,i) \longmapsto  {H}^{i}(a_{!}\,u_{\ast}\,u^{\ast}\,a^{!}(\mathcal O_X))
\]
where now $a_{!}, u_{\ast}, u^{\ast}, a^{!}$ are the operations in the derived category of the category of complexes of $\mathcal D$-modules with regular holonomic cohomology and $H$ is ordinary cohomology. For this definition to make sense, one has to consider $\mathcal D$-modules on singular varieties, and also pull-backs and push-forwards for them, see digression 2. We denote $\mathfrak N\mathfrak I (X)$ the corresponding diagram category and $\mathcal F_{\mathcal D}: \mathfrak{NI}(X) \longrightarrow rh\mathcal D_{X}-\mathfrak{mod}$ the functor obtained via Nori-Ivorra's theorem. An object in $rh\mathcal D_{X}-\mathfrak{mod}$ will be said to be motivic if it is in the essential image of $\mathcal F_{\mathcal D}$, similarly for morphisms.
\end{enumerate}
\end{Remarks}

It will be convenient to have some more understanding of the functor $\mathcal F_{\mathcal D}: \mathfrak{NI}(X) \longrightarrow rh\mathcal D_{X}-\mathfrak{mod}$. Since the category $rh\mathcal D_X-\mathfrak{mod}$ is finite and Hom-finite, by Ivorra's theorem there is an equivalence between this category and the category of comodules over a $\mathbb C$-coalgebra which are finite dimensional over $\mathbb C$. Then we have an exact, faithful functor $\omega: rh\mathcal D_X-\mathfrak{mod} \longrightarrow 
\mathbb C-\mathfrak{fvs}$, the composition of this equivalence and the forgetful functor. 

Let $D(rh\mathcal D_X-\mathfrak{mod})$ be the diagram
of the category $rh\mathcal D_X-\mathfrak{mod}$. Then, $\omega$ can be viewed as a representation of $D(rh\mathcal D_X-\mathfrak{mod})$ in $\mathbb C-\mathfrak{fvs}$, and by Nori's theorem (in the version in \cite{Ivo}*{Theorem 3.1}), $rh\mathcal D_X-\mathfrak{mod}$ is equivalent to the diagram category $\mathcal C(D(rh\mathcal D_X-\mathfrak{mod}), \omega)$, which in turn is equivalent to the category $A(D(rh\mathcal D_X-\mathfrak{mod}), \omega)-\mathfrak{comod}$.

It follows from the proof of \cite{Ivo}*{Lemma 5.1} that we can take $\mathfrak N\mathfrak I (X)$ 
equal to the category $\mathcal C(\mathfrak D^{Nori}_X, \omega \circ T_{\mathcal D})$, and $\widetilde{T_{\mathcal D}}: \mathfrak D^{Nori}_X \longrightarrow  \mathcal C(\mathfrak D^{Nori}_X, T_{\mathcal D})$ equal to the functor $\widetilde{\omega\circ T_{\mathcal D}}: \mathfrak D^{Nori}_X \longrightarrow \mathcal C(\mathfrak D^{Nori}_X, \omega \circ T_{\mathcal D})$. The situation is summarized in the following diagram 
\[
\xymatrix{ {} & 
{\ \ \ \mathfrak{NI}(X)=\mathcal C(\mathfrak D^{Nori}_X, \omega \circ T_{\mathcal D})} \ar@/_0.1 pc/[1,1]^-{\mathcal F_{\mathcal D}}
\ar@/^0.1pc/[2,0]^-(.35){\mathcal F_{\omega\circ T_{\mathcal D}}=for}
& {} \\ {\mathfrak{D}^{Nori}_X} \ar@/^0pc/[-1,1]^-{\widetilde{T_{\mathcal D}}\,=\,\widetilde{\omega\circ T_{\mathcal D}\ \ }}
\ar@/^0.1pc/[0,2]_-{\!\!\!\!\!\!\!\!\!\!\!\!\!\!\!\!\!\!\!\!\!\!\!\!\!\!\!\!     T_{\mathcal D}}
\ar@/_0.1pc/[1,1]_-{\omega\circ T_{\mathcal D}}
& {} 
& {rh\mathcal D_{X}-\mathfrak{mod}
.} \ar@/_0pc/[1,-1]^-{\omega}\\ {} 
& {\mathbb C-\mathfrak{fvs} } 
   & {}}
\]
By loc. cit., the functor $\mathcal F_{\mathcal D}$ is characterized by the existence of natural equivalences 
(invertible $2$-morphisms) 
\[
 \alpha: \mathcal F_{\mathcal D}\circ \widetilde{T_{\mathcal D}} \stackrel{\sim}\longrightarrow T_{\mathcal D}  \text{\ \ and\ \  } \delta: \omega\circ \mathcal F_{\mathcal D} \stackrel{\sim}\longrightarrow \mathcal F_{\omega\circ T_{\mathcal D}}
\]
veryfing that the square of functors in \cite{Ivo}*{pg. 1112} is commutative. If we apply the functoriality explained in \cite{Ivo}*{3.3}
in this case (in the notation of loc. cit., take $\mathcal Q_1=\mathfrak{D}^{Nori}_X$ and $\mathcal Q_2= D(rh\mathcal D_X-\mathfrak{mod})$, $T_1=\omega\circ T_{\mathcal D}$, $T_2=\omega$) then we obtain a morphism of coalgebras
\[
 A(\mathfrak{D}^{Nori}_X, \omega\circ T_{\mathcal D})\longrightarrow  A(D(rh\mathcal D_X-\mathfrak{mod}), \omega),
\]
and the corresponding functor of (co)restriction of scalars 
\begin{eqnarray*}
\mathfrak{NI}(X)\cong A(\mathfrak{D}^{Nori}_X, \omega\circ T_{\mathcal D})-\mathfrak{comod} \longrightarrow \\ A(D(rh\mathcal D_X-\mathfrak{mod}), \omega)-\mathfrak{comod}\cong rh\mathcal D_X-\mathfrak{mod}
\end{eqnarray*}
and one can easily check that, up to isomorphism, this is the functor $\mathcal F_{\mathcal D}$ we wanted to understand.

{\bf 3. The relation with local cohomology} 

In this section we will consider mainly affine varieties, in this case $\mathcal D$- and $\mathcal O$-modules can be identified to their modules of global sections, what we will do without explicit mention.
Let $R$ be the coordinate ring of a complex affine variety, let $I\subset R$ be an ideal. Put $X=\text{Spec}\,R$, denote $Y\subset X$ the zero-locus of $I$. It is known that the local algebraic cohomology module $H^i_Y(\mathcal O_X)$ (alias $H^i_I(R)$)
is a regular holonomic $\mathcal D_X$-module for all $i\geqslant 0$. 
\newpage


\begin{Lemma}
  $H^i_Y(\mathcal O_X)$ is motivic \footnote{Or, if one prefers the language of perverse sheaves, $^{p}H^i(R\Gamma_Y(\mathbb C_X))$ is motivic.}.
\end{Lemma}

\begin{Proof}
 The vertex of $\mathfrak D^{Nori}_X$ corresponding to $(Y,\emptyset ,i)$ maps to $H^i_Y(\mathcal O_X)$  by $T_{\mathcal D}$ (see e.g. \cite{Bo}*{VI, Theorem 7.13}).
 By Nori's theorem, $T_{\mathcal D}= F_{\mathcal D} \circ \widetilde{T}$, thus $\widetilde{T}(Y,\emptyset ,i)$ is the motive of $H^i_Y(\mathcal O_X)$.
\end{Proof}
\medskip

\begin{Remarks}

\begin{enumerate}
\item The functor $\mathcal F_{\mathcal D}$ is faithful but we do not know if it is fully faithful. Thus, if we have a regular holonomic $\mathcal D_X$-module $M$ which is in the essential image of $\mathcal F_{\mathcal D}$, we cannot say which is ``its motive'', because a priori there could exist non-isomorphic objects in
$\mathfrak{NI}(X)$ which map to $M$. However, for the local cohomology modules considered above, notice that the given argument provides a distinguished preimage. 
\item What is ``the motive of $H^i_I(R)$''? As a set, it is $H^i_I(R)$, but we regard it not only as a $\mathcal D_X$-module (or as a $R$-module) but as a $A(\mathfrak{D}^{Nori}_X, \omega\circ T_{\mathcal D})$-comodule (or as a $\text{End}_{\mathbb C}(\omega \circ T_{\mathcal D})$-module, see digression 1).
 
\item 
Let $Z\subset Y$ be two closed subvarieties
of $X$, denote $i:Y\hookrightarrow X$, $j: Z \hookrightarrow X$ the inclusions. The isomorphism $\varphi_Y: i_{\ast}i^{!}(\mathcal O_X) \longrightarrow \mathbb R \Gamma_Y(\mathcal O_X)$ (in the derived category of 
the category of $\mathcal D_X$ modules with regular holonomic cohomology, see \cite{Bo}*{VI, Theorem 7.13}) is functorial, in the sense that we have a commutative square
 \[ \begin{CD}
j_{\ast}j^{!}(\mathcal O_X) @>>> i_{\ast}i^{!}(\mathcal O_X)\\
@V\varphi_{Z}VV @VV\varphi_{Y}V\\
\mathbb R \Gamma_{Z}(\mathcal O_X)@>>> \mathbb R \Gamma_{Y}(\mathcal O_X)\,.
\end{CD} \]

Here, $\mathbb R \Gamma_{Z}(\cdot)$ is the derived functor of the functor on sheaves of abelian groups given by $\mathcal F \longmapsto \Gamma_{Z}(\mathcal F)$, where for $U\subset X$ open, $\Gamma_{Z}(\mathcal F)$ is given by $U\longmapsto \Gamma_{Z\cap U}(U,\mathcal F_{\mid U})$,  and similarly for $Y$. Denote $\mathbb R \Gamma_{Y/Z}(\mathcal O_X)$ the derived functor of the functor 
$\mathcal F \longmapsto \Gamma_{Y/Z}(\mathcal F)$, where $\Gamma_{Y/Z}(\mathcal F)$ is the sheaf defined by $U\longmapsto \Gamma_{U\cap Y}(U,\mathcal F_{\mid U})/\Gamma_{U\cap Z}(U,\mathcal F_{\mid U})$, considered in \cite{RD}*{Chapter IV}. We have distinguished triangles
\begin{eqnarray}
 \mathbb R \Gamma_{Z}(\mathcal O_X) \longrightarrow \mathbb R \Gamma_{Y}(\mathcal O_X) \longrightarrow \mathbb R \Gamma_{Y/Z}(\mathcal O_X) \stackrel{+1}\longrightarrow\\
 j_{\ast}j^{!}(\mathcal O_X) \longrightarrow  i_{\ast}i^{!}(\mathcal O_X) \longrightarrow i_{\ast}\,u_{\ast}\,u^{!}\,i^{!}(\mathcal O_X) \stackrel{+1}\longrightarrow \,,
\end{eqnarray}
where $u:Y\smallsetminus Z \hookrightarrow Y$ is the inclusion. Triangle (1) is tautological, for (2) see \cite{Ivo2}*{Lemma 3.4}.
Then, it follows from the functoriality of the isomorphisms $\varphi_Y$ that, for $\ell\geqslant 0$, we have isomorphisms of $\mathcal D_X$-modules
\[
T_{\mathcal D}(Y,Z,\ell)={H}^{\ell}(i_{\ast}\,u_{\ast}\,u^{!}\,i^{!}(\mathcal O_X))\cong H^{\ell}(\mathbb R \Gamma_{Y/Z}(\mathcal O_X)) = H^{\ell}_{Y/Z}(\mathcal O_X) \,,
\]
where the notation on the right is taken from \cite{RD}. Thus, $\widetilde{T}_{\mathcal D}(Y,Z,\ell)$ is ``the motive of $H^{\ell}_{Y/Z}(\mathcal O_X)$''.
%
%
\end{enumerate}
\end{Remarks}
\medskip

\begin{Example}
 A baby motivic theorem\footnote{See more grown-up examples in \cite{Ivo2}, see also section 4 below.}: Let $I,J\subset R$ be ideals. Consider the Mayer-Vietoris exact sequence
 \[
  \dots \longrightarrow  H^i_{I+J}(R) \longrightarrow H^i_I(R)\oplus H^i_J(R)\longrightarrow H^i_{I\cap J}(R)\longrightarrow H^{i+1}_{I+J}(R) \longrightarrow \dots
 \]
We claim that this sequence can be promoted to an exact sequence in $\mathfrak{NI}(X)$: The terms in the sequence are motivic as we have just seen, the
morphisms are also motivic: Let $Z_{I+J}\hookrightarrow Z_{I\cap J}$ be the subvarieties of $X$ defined by $I+J,I\cap J$. Then, the connecting homomorphisms in the sequence above come from the composition of the morphisms corresponding to the edges $(Z_{I\cap J},\emptyset,i)\longrightarrow (Z_{I\cap J}, Z_{I+J}, i) \longrightarrow
(Z_{I+J}, \emptyset, i+1)$ in the diagram $\mathfrak D^{Nori}_X$. The other morphisms come 
from arrows of type i). Since the functor $\mathcal F_{\mathcal D}: \mathfrak{NI}(X) \longrightarrow \mathcal D_X-\mathfrak{mods}$ is faithful,
the corresponding long sequence in $\mathfrak{NI}(X)$ will be exact. $\Box$
\end{Example}

Assume $R=\mathbb C[x_1,\dots,x_n]$ for the rest of this section, put $\mathfrak m=\langle x_1,\dots,x_n\rangle \subset R$, denote $\mathbf 0$ the origin in $X=\text{Spec R}$. In \cite{Lyu}, G. Lyubeznik proved that the Bass numbers $\lambda_{r,i}:=\mu_r(\mathfrak m, H^{n-i}_I(R_{\mathfrak m}))$ are finite and depend only on the ring $\mathcal O_{Y,\mathbf{0}}=R_{\mathfrak m}/I_{\mathfrak m}$, since then they have been called Lyubeznik numbers in the literature. One has
\[
 H_{\mathbf 0}^{r}(H_Y^{n-i}(\mathcal O_X))\cong H_{\bf 0}^n(\mathcal O_X)\,\oplus\, \stackrel{\lambda_{r,i})}\dots \, \oplus \,H_{\bf 0}^n(\mathcal O_X).
\]
The $\mathcal D_X$-modules $H_{\bf 0}^n(\mathcal O_X)$ are simple, so $\lambda_{r,i}$ is the length of $ H_{\mathbf 0}^{r}(H_Y^{n-i}(\mathcal O_X))$ as a $\mathcal D_X$-module.
\medskip

\begin{Proposition}
The $\mathcal D_X$-modules $H^{r}_{\mathbf 0}(H_Y^{\ell}(\mathcal O_X))$ are motivic for all $r,\ell\in\mathbb N$.
\end{Proposition}
\begin{Proof} 
For $J=\{j_1,\dots,j_k\}\subset\{1,\dots,n\}$ a multi-index with $j_1<...<j_k$, put $\lvert J\rvert=k$ and denote by $V_J$ the open subset of $X$
defined by $x_{j_1}\dots x_{j_k}\neq 0$, let $v_J:V_J \hookrightarrow X$ be the inclusion map.

The local algebraic cohomology $H_{\mathbf 0}^{\ast}(H_Y^\ell(\mathcal O_X))$ is the cohomology of the Cech complex
\[
\cdots \longrightarrow \oplus_{\mid J \mid =k} (v_J)_{\ast}v_J^{\ast}(H_Y^\ell(\mathcal O_X)) \longrightarrow \oplus_{\mid J \mid =k+1} (v_J)_{\ast}v_J^{\ast}(H_Y^\ell(\mathcal O_X)) \longrightarrow\dots 
\]
Since the category of perverse Nori motives is abelian, it is enough to show that both the terms and the morphisms in the above complex are motivic. 
Let $U_J=V_J\cap Y$, $u_J: U_J \hookrightarrow Y$, $a:Y\hookrightarrow X$ the inclusion maps, $i_J=a\circ u_J: U_J \longrightarrow X$.

The cohomology of the complex $a_{!}(u_J)_{\ast}(u_J)^{\ast}a^{!}(\mathcal O_X)$ is motivic just by definition of $T_{\mathcal D}$. We have
\[
 (i_J)_{\ast}(i_J)^{\ast}a_{!}\,a^{!}(\mathcal O_X)=a_{\ast}(u_J)_{\ast}(u_J)^{\ast}a^{\ast}\,a_{!}\,a^{!}(\mathcal O_X).
\]
Since $a$ is a closed immersion, $a^{\ast}a_{!}=a^{\ast}a_{\ast}=id$, and it follows that \\ 
$(i_J)_{\ast}(i_J)^{\ast}a_{!}\,a^{!}(\mathcal O_X)$ is motivic.
Since $a_{!}\,a^{!}(\mathcal O_X)$ is supported on $Y$, 
\[
  (i_J)_{\ast}(i_J)^{\ast}a_{!}\,a^{!}(\mathcal O_X)= (v_J)_{\ast}(v_J)^{\ast}a_{!}\,a^{!}(\mathcal O_X).
\]
Also, $a_{!}\,a^{!}(\mathcal O_X)=a_{\ast}\,a^{!}(\mathcal O_X)=R\Gamma_Y(\mathcal O_X)$, 
and $(v_J)_{\ast}, (v_J)^{\ast}$ are derived functors of exact functors because $v_J$ is an affine open embedding. It follows that for the terms in the Cech complex
above we have
\[
(v_J)_{\ast}(v_J)^{\ast}H^\ell(R\Gamma_Y(\mathcal O_X)) = H^\ell((u_J)_{\ast}(u_J)^{\ast}a_{!}\,a^{!}(\mathcal O_X))\,,
\]
and so they are images of vertices of $\mathfrak{D}_{Nori}$ by $T_{\mathcal D}$, 
thus they are motivic.

Let $J=\{j_1,\dots,j_k\}\subset\{1,\dots,n\}$ be a multi-index as above, put $J'=J\cup\{j\}$ with $j\not\in J$. 
Then in $\mathfrak D^{Nori}_X$ there are edges  $(Y,Y-U_J,i) \longrightarrow (Y,Y-U_{J'},i)$ for all $i\geqslant 0$, and 
the morphisms in the Cech complex are linear combinations of morphisms corresponding to these edges via Ivorra's representation.
$\Box$
\end{Proof}
\medskip

\begin{Remarks}
\begin{enumerate}
\item Notice that the argument above gives also a distinguished object of $\mathfrak{NI}(X)$ mapping to $H^{r}_{\mathbf 0}(H_Y^{n-i}(\mathcal O_X))$ via the faithful, exact functor $\mathcal F_{\mathcal D}$, denote it 
$[H^{r}_{\mathbf 0}(H_Y^{n-i}(\mathcal O_X))]$.
Since $H^{r}_{\mathbf 0}(H_Y^{n-i}(\mathcal O_X))$ is a $\mathcal D_X$-module of finite length $\lambda_{r,i}$, $[H^{r}_{\mathbf 0}(H_Y^{n-i}(\mathcal O_X))]$ is also an object of finite length in $\mathfrak{NI}(X)$, and we have a bound
\[
\text{length}([H^{r}_{\mathbf 0}(H_Y^{n-i}(\mathcal O_X))])\leqslant \lambda_{r,i}\,.
\]
I do not know if this is an equality except in some very simple cases, for example if $\lambda_{r,i}=1$. 
The numbers in the left hand side of the inequality above might be called ``motivic Lyubeznik numbers''. I do not know if they depend only on the local ring $\mathcal O_{Y,\mathbf 0}$.

\item In proposition 1 the initial datum is an ideal on a polynomial ring. Since not all singularity germs can be globalized, it would be desirable
 to have a proof for ideals in a power series ring with complex coefficients\footnote{A possibility would be to work with $\mathcal D$-modules on formal schemes. Beilinson and Drinfeld introduced in \cite{QHE} a definition of $\mathcal D$-modules even on ind-schemes, with direct and inverse images, so it is likely that their work gives a solution to the question posed above.}. 
\end{enumerate}
\end{Remarks}

Notice that we could have taken any subdiagram of $\mathfrak D^{Nori}_X$ to define a diagram category, in particular we could take only those vertices 
corresponding to triples $(Y\hookrightarrow X,Z,i)$ where $Y\hookrightarrow X$ is a closed immersion, denote this diagram $\mathbf{Lcoh}_X$. Such a choice would give a diagram category, say $\mathfrak{Lcoh}_X$, and a functor $\mathcal F_{\mathbf{Lcoh}_X}:\mathfrak{Lcoh}_X \longrightarrow rh\mathcal D_{X}-\mathfrak{mod}$. They might be  appropriate for the study of local cohomology modules.
In particular, the category $\mathfrak{Lcoh}_X$ allows to make precise the discussion in the introduction to these notes: It seems reasonable to interpret a ``functorial, additional structure on local cohomology modules''\footnote{Additional means here additional to the structure of $\mathcal D_X$-module.} as a $\mathbb C$-linear, abelian category $\mathcal A$, endowed with  a representation $\mathbf{Lcoh}_X \longrightarrow \mathcal A$ and an exact, faithful functor $\mathcal  A\longrightarrow rh\mathcal D_{X}-\mathfrak{mod}$\footnote{This functor would correspond to forget the additional structure and retain only the one of $\mathcal D_X$-module.}, such that the triangle of functors 
\[
    \begin{tikzcd}
     {} &  \mathcal A \arrow{dr} \\
    \mathbf{Lcoh}_X \arrow{ur} \arrow{rr}{T_{\mathcal D\mid \mathbf{Lcoh}_X }} && rh\mathcal D_{X}-\mathfrak{mod}
    \end{tikzcd}
\] 
commutes. Assuming we are in such a situation, by the theorem of Nori-Ivorra, there is a faithful, exact functor $\mathfrak{Lcoh}_X \longrightarrow \mathcal A$ making the corresponding diagram of functors commutative. So, roughly, ``the richest'' functorial, additional structure that we can put on the local cohomology modules $H^i_Y(\mathcal O_X)$ would be that of a $A(\mathbf{Lcoh}_X, T_{\mathcal D\mid \mathbf{Lcoh}_X })$-comodule. 

It seems a natural question to ask ``how far'' is the diagram category $\mathfrak{Lcoh}_X$ from its essential image in the category $rh\mathcal D_{X}-\mathfrak{mod}$ (or, equivalently,  how far is the functor $\mathcal F_{\mathbf{Lcoh}_X}$
from being full). 
I do not know the answer, which might depend on the variety $X$. 
\medskip

{\bf 4. Basic lemma. Stratifications.}
 
 In \cite{Ivo2}*{Lemma 4.5}, Ivorra uses the so-called ``basic lemma'' of Beilinson to prove a lemma, a special case of which is the following proposition:
 
 \begin{Proposition}
 Let $X$ be a smooth\footnote{Most likely, smoothness is not needed, but I have not fully checked it.} affine complex manifold of dimension $d$, $Y \subset X$ a subvariety of dimension $n$. Given 
 a closed subvariety $W\subset Y$ such that $\dim(W)\leqslant n-1$ there exists a closed subvariety $Z\subset Y$ such that $W\subset Z$, $\dim(Z)\leqslant n-1$ and for every $j\neq d-n$ one has\footnote{In the notations of \cite{Ivo2}, one has $TH^{\mathcal M}_X(Y,Z,i)=0$ for $i\neq n$. The shifting on the vanishing range is due to the fact that, if $\mathcal M$ is the category of perverse sheaves, then $T_{\mathcal D}(Y,Z,j)$ corresponds to $TH^{\mathcal M}_X(Y,Z,d-j)$ via the Riemann-Hilbert correspondence.}
 \[
  H^j_{Y/Z}(\mathcal O_X)=0.
 \]
 If $Y$ is reduced and irreducible, one can choose $Z$ such that its open complement is smooth.
\end{Proposition}

From the long exact sequence coming from the tautological triangle (1) in Remark iii) above, one gets, in algebraic terms:

\begin{Proposition}
Let $R$ be the coordinate ring of a complex affine manifold of dimension $d$, $I\subset R$ an ideal of height $h$. Let $J\supset I$ be an ideal such that
$ht(J)\geqslant h+1$. Then, there is an ideal $\mathfrak a$ with $I\subset \mathfrak a \subset J$ and $ht(\mathfrak a) \geqslant h+1$, such that one has:
\begin{enumerate}
 \item An isomorphism of $R$-modules $H^i_{I}(R)\cong H^i_{\mathfrak a}(R) $ for every $i\neq h, h+1$.
\item An exact sequence
\begin{eqnarray*}
 0 \longrightarrow H^h_{\mathfrak a}(R)\longrightarrow  H^{h}_{I}(R) &\longrightarrow& 
 H^h_{\mathfrak a / I}(R)
 \longrightarrow \\  H^{h+1}_{\mathfrak a}(R)&\longrightarrow&  H^{h+1}_{I}(R) \longrightarrow 0.
\end{eqnarray*}
\end{enumerate}
\end{Proposition}

While Beilinson's basic lemma is quite a deep result, the above propositions are derived from rather special cases of it, so it is possible that one can prove them directly
in a simpler way (and probably they hold also in positive characteristic). In any case, the motivic point of view allows to deduce them automatically from the result of Beilinson and Ivorra. 

Proposition 1 is used in \cite{Ivo2} to construct cellular stratifications, we rephrase the construction in the present setting. In the sequel, all varieties are assumed to be closed subvarieties of a fixed smooth algebraic variety $X$.

\begin{Definition}
Let $Y$ be a complex algebraic variety of dimension $d$. A {\it stratification} $Y_{\bullet}$ of $Y$ is a sequence of closed subsets
\[
 Y_{\bullet}: \emptyset=Y_{-1}\subset Y_0\subset \dots Y_i\subset Y_{i+1}\subset \dots \ \ \ i\in\mathbb N\cup\{-1\}
\]
such that $\dim Y_i\leqslant i$ and $Y_m=Y$ for some $m\in\mathbb N$.

The stratification is {\it cellular} if the following two conditions are fulfilled:
\begin{enumerate}
 \item If $\dim Y_i=i$, then for every $k\neq d-i$ one has $H^k_{Y_i/Y_{i-1}}(\mathcal O_X)=0$.
 \item If $\dim Y_i<i$, then $Y_i=Y_{i-1}$ (and then $H^k_{Y_i/Y_{i-1}}(\mathcal O_X)=0$ for all $k\in\mathbb N$). 
\end{enumerate}
\end{Definition}
\medskip

\begin{Proposition}(\cite{Ivo2}*{Corollary 4.8})
 There always exist cellular stratifications of $Y$. 
\end{Proposition}

If we have a topological space with a filtration by sheaves of families of supports, the cohomology of a complex of abelian sheaves can be computed by means of a spectral sequence (see \cite{RD}). If the filtration is given by a cellular stratification, then:
\medskip

\begin{Proposition}(A special case of \cite{Ivo2}*{Proposition 4.12})
 Let $Y_{\bullet}$ be a cellular stratification of $Y$. Then, there is a complex $\mathfrak S_{\ast}(Y_{\bullet})$
\[
 \longrightarrow  H^{d-i}_{Y_i/Y_{i-1}}(\mathcal O_X) \longrightarrow  H^{d-i+1}_{Y_{i-1}/Y_{i-2}}(\mathcal O_X)\longrightarrow \dots\longrightarrow  H^{d}_{Y_0/Y_{-1}}(\mathcal O_X)\longrightarrow 0\,,
\]
 where $H^{0}_{Y_d/Y_{d-1}}(\mathcal O_X)$ is placed in degree $0$, such that 
 \[
  H^k_Y(\mathcal O_X) \cong H^k(\mathfrak S_{\ast}(Y_{\bullet})).  
 \]
\end{Proposition}
The morphisms in the complex $\mathfrak S_{\ast}(Y_{\bullet})$ are the following ones:
There is a tautological distinguished triangle
\[
\mathbb R \Gamma_{Y_{i-1}/Y_{i-2}}(\mathcal O_X) \longrightarrow \mathbb R \Gamma_{Y_i/Y_{i-2}}(\mathcal O_X) \longrightarrow \mathbb R \Gamma_{Y_i/Y_{i-1}}(\mathcal O_X) \stackrel{+1}\longrightarrow
\]
which gives rise to a corresponding long exact sequence
\[
 \dots \longrightarrow H_{Y_i/Y_{i-2}}^l(\mathcal O_X) \longrightarrow H_{Y_i/Y_{i-1}}^l(\mathcal O_X) \longrightarrow H_{Y_{i-1}/Y_{i-2}}^{l+1}(\mathcal O_X) \longrightarrow \dots
\]
Since the stratification is cellular, there is only one non-zero connecting homomorphism $H_{Y_i/Y_{i-1}}^{d-i}(\mathcal O_X) \longrightarrow H_{Y_{i-1}/Y_{i-2}}^{d-i+1}(\mathcal O_X)$, and this is the $(d-i)$-th differential in $\mathfrak S_{\ast}(Y_{\bullet})$. The complex in Proposition 5 resembles much the Cousin complexes in \cite{RD}.
\medskip

{\bf 5. Positive characteristic} 

We assume now that $k$ is a perfect field of characteristic $p>0$. We will use some non-explained terminology, we refer for it to \cite{EK} and \cite{Oh}.

Let $X$ be a separated $k$-scheme of finite type over $k$, assume there is a proper smooth $k$-schema $P$ and an immersion $X\hookrightarrow P$. M. Emerton and M. Kisin introduced in \cite{EK} the bounded derived category $D^b(\mathcal D_{P,F})$ of the category of $\mathcal D_P$-modules with Frobenius structure, and the subcategory $D^b_{lfgu}(\mathcal D_{P,F})^{\circ}$ of $D^b(\mathcal D_{P,F})$ consisting of complexes with locally finitely generated unit cohomology and finite $Tor$-dimension over $\mathcal O_P$.

In his thesis, S. Ohkawa proved that
the subcategory of $D^b_{lfgu}(\mathcal D_{P,F})^{\circ}$ of complexes supported on $X$ does not depend on $P$\footnote{Up to natural equivalence. At this point one could also use the method of Saito \cite{Sai} to obtain a single category and not just a category up to equivalence, see digression 2.}, denote this subcategory by $D^b_{lfgu}(X)$. He also established the existence of direct images, extraordinary inverse images and tensor products (see \cite{Oh}, the setting and results of Ohkawa are in fact more general). The Riemann-Hilbert correspondence proved in \cite{EK} extends to this setting, that is, there is an anti-equivalence between $D^b_{lfgu}(X)$ and $D(X_{\text{\'et}},\mathbb Z/p\mathbb Z)$, given by an analogue of the solutions functor. There is no analogue of the De Rham functor, because duality fails. It follows from this version of the Riemann-Hilbert correspondence that, in the sequel, we could have worked with constructible $\mathbb Z/p\mathbb Z$-sheaves on the \'etale site of $X$.

Both Emerton-Kisin and Ohkawa work in the more general setting of schemes over rings $W_n(k)$ of truncated Witt vectors. We will only consider the case $n=1$, in this case there is an equivalence between the category of unit $\mathcal D_{P,F}$-modules and the category of unit $\mathcal O_{P,F}$-modules, preserving the condition of being locally finitely generated. Moreover, there is also an equivalence for the derived categories of the categories of complexes with locally finite generated unit cohomology which is compatible with cohomological operations and support conditions (see \cite{EK}*{\S 15.2, 15.4.3}).  

The standard t-structure in $D^b_{lfgu}(X)$ has as heart the category $\mu_{lfgu}(X)$ consisting of locally finitely generated unit $\mathcal D_{P,F}$-modules supported at $X$\footnote{They might be regarded as an analogue of the regular holonomic modules in the complex case.}. This category is finite (follows from \cite{Lyu2}*{Theorem 3.2}) and Hom-finite (this follows from \cite{Ho}*{Theorem 5.1}). To prove these two claims one has to keep in mind that, in the affine case, locally finitely generated unit $\mathcal O_{F,X}$-modules are the same as the F-finite modules introduced by Lyubeznik in \cite{Lyu} .

This suggest that we can partially reproduce Ivorra's construction in positive characteristic. Let $X$ be a quasi-projective $k$-manifold and consider the diagram $\mathfrak D'_{Nori}$
where the vertices are triples $(Y \stackrel{a}\to X,Z,i)$, $Z$ closed in $Y$, but this time we impose also that $Y$ is quasi-projective and $a:Y\longrightarrow X$ is proper.
This is because in the setting of Emerton-Kisin and Ohkawa, for a general morphism $f$ only $\mathbb Rf_{\ast}$ and $f^{!}$ are defined.
The set of edges is the same as in $\mathfrak D^{Nori}_X$, except that in the edges of type i) we require $f$ to be a proper map. Then we consider the assignment $T'_{\mathcal D}: \mathfrak D'_{Nori} \longrightarrow 
D^b_{lfgu}(X)$ given by\footnote{As in the complex setting, if $f$ is a morphism of schemes we will denote by $f_{\ast}$ the functor $\mathbb Rf_{\ast}$. Notice that, 
what makes this definition the analogue of Ivorra's in the complex case, is the fact that for constructible sheaves we have $f_{!}=f_{\ast}$ if $f$ is proper and $f^{*}=f^{!}$ if $f$ is \'etale (in particular for an open embedding).}
\[
T'_{\mathcal D}: (Y,Z,i) \longmapsto  {H}^{i}((a_{\mid Y \setminus Z})_{\ast}\,(a_{\mid Y\setminus Z})^{!}(\mathcal O_X))\,.
\]
Looking at Ivorra's proof in \cite{Ivo} we observe that, what we need to have a true representation of the diagram  $\mathfrak D'_{Nori}$ is to prove the following statements:
\begin{enumerate}

 \item (Adjunction morphisms):
 
 Given $f:Y \longrightarrow Z$, there is a natural transformation $f_{\ast}f^{!}\longrightarrow id$. If $k:U\hookrightarrow Z$ is an open immersion, then  this transformation is an isomorphism of functors. 
 
 \item (Connection with local cohomology): If $i: Y\hookrightarrow Z$ is a closed immersion of smooth quasi-projective $k$-varieties and $M^{\bullet}$ is an object in $D^b_{lfgu}(Z)$, then
\[
 i_{\ast}i^{!}M^{\bullet}=\mathbb R\Gamma_Y(M^{\bullet}).
\]
 
 \item (Localization triangle) If $Y\subset Z$ is closed, denoting $i: Y\hookrightarrow Z$ and $j:Z\smallsetminus Y \hookrightarrow Z$ the inclusions, for each object $M$ of $D^b_{lfgu}(Z)$ there is an exact triangle
 \[
  i_{\ast}i^{!}M\longrightarrow M \longrightarrow j_{\ast}j^{!}M \stackrel{+1}\longrightarrow
 \]
 
 \item (Base change morphism): Given a Cartesian square 
 \begin{equation*}
  \xymatrix@R+2em@C+2em{
  V \ar[r]^-l \ar[d]_-h & Y \ar[d]^-f \\
  Z \ar[r]_-g & T
  }
 \end{equation*}
there is a natural transformation $l_{\ast}\,h^{!}\longrightarrow f^{!}\,g_{\ast}$.

\end{enumerate}

{\it Proof of i)-iv)}: We use the definitions and the notations of \cite{Oh}. 
%
%
The first adjunction in i) follows from the definitions and the smooth case (\cite{EK}*{Theorem 4.4.1}).
If $Z$ is smooth, then the second claim in i) follows from \cite{EK}*{Lemma 4.3.1}. In general, assume  $k:U\hookrightarrow Z$ is an open immersion. Choose $P$ proper and smooth and $V\subset P$ open such that there is a closed inmersion $Z \hookrightarrow V$. Let $W$ be open in $U$ such that $U=W\cap Z$.
We have a commutative square
 
  \begin{tikzcd}
    & U\arrow[hookrightarrow]{r}\arrow[hookrightarrow]{d}{k} 
   & W\arrow[hookrightarrow]{r}{j}
    & P\arrow{d}{id}
    \\
    & Z\arrow[hookrightarrow]{r}
    & U\arrow[hookrightarrow]{r}
    & P
  \end{tikzcd}
  
and  $k^{!}=j_{\ast}\mathbb R\Gamma_U j^{!}:\mathcal C_{P,Z} \longrightarrow \mathcal C_{P,W}$ by Ohkawa's definition of $k^{!}$.
Then, the claimed isomorphism follows from the same one in the smooth case applied to $j$. 

For ii), choose again $P$ proper and smooth and $U\subset P$ open such that there is a closed inmersion $Z \hookrightarrow U$. Then we have

  \begin{tikzcd}
    & Y\arrow[hookrightarrow]{r}\arrow[hookrightarrow]{d}{i} 
    & U\arrow[hookrightarrow]{r}{j}
    & P\arrow{d}{id}
    \\
    & Z\arrow[hookrightarrow]{r}
    & U\arrow[hookrightarrow]{r}
    & P\ .
  \end{tikzcd}
 
Let $D\subset P$ closed such that $Y=U\cap D$. By \cite{Oh}*{Proposition 3.5} we have $\mathbb R \Gamma_Y j^{!}=j^{!}\mathbb R\Gamma_D$, thus given an object $M$ in $\mathcal C_{P,Z}$ we have
\[
 i_{\ast}i^{!}M=j_{\ast}\mathbb R \Gamma_Y j^{!}(M)=j_{\ast}j^{!}\mathbb R\Gamma_D(M)=\mathbb R\Gamma_D(M),
\]
where the last isomorphism follows form i). But $\mathbb R\Gamma_{P\smallsetminus U}(M)=0$, so \\ $\mathbb R\Gamma_D(M)=\mathbb R\Gamma_Z(M)$.

Item iii) follows from ii). Item iv) is well-known for $\mathcal D$-modules, see e.g. \cite{HKT}*{Proposition 1.7.3}\footnote{In \cite{HKT} it is proved that the base change morphism is in fact an isomorphism, but this is not necessary for our purposes. Other proofs of similar base changes use the full six-functor formalism and so cannot just be transcribed into the positive characteristic setting, for example the one given in \cite{KS}*{Chapter III, Proposition 3.1.9} for constructible sheaves.}, in the present setting the result can be proved following the same steps as in loc. cit., we recall a few details: First, decomposing $f$ via the graph embedding, we can assume $f$ is either a closed immersion or a projection. If $f$ is a closed immersion, then there is an equivalence $f^{!}f_{\ast}\cong id_Z$ (in the smooth case this is the analog of Kashiwara's equivalence, which is proved in \cite{EK}*{Proposition 15.5.3}, in Ohkawa's setting this is a consequence of the definitions, because Kashiwara's equivalence is built into them). Then we have
natural transformations
\[
 l_{\ast}\, h^{!} \simeq f^{!}\,f_{\ast}\,l_{\ast}\,h^{!}=f^{!}\,g_{\ast}\,h_{\ast}\,h^{!}  \rightarrow f^{!}\,g_{\ast}
\]
where the last one comes from i) applied to $h$. The case of a projection is proved as in \cite{HKT}.
%
%
%
%
%
$\Box$

Then, as in the complex case, one can consider the corresponding category of motives $\mathfrak{NI}(X)$. If $R=k[x_1,\dots,x_n]$,
$I\subset R$ is an ideal and $\ell\geqslant 0$, we obtain that the modules $H^\ell_I(R)$ are motivic and also that ``motivic Lyubeznik numbers'' can be defined\footnote{Notice that,
in the proof of Proposition 1 the morphisms in the Cech complex correspond to edges given by proper maps, in fact identities. Thus, the iterated local cohomologies in loc. cit. are also motivic in positive characteristic.}. As before, they are bounded above by the usual Lyubeznik numbers.

\begin{Remark}
It is likely that a similar formalism of Nori's categories, etc. can be established for Blickle and B\"ockle Cartier crystals over $F$-finite schemes (see \cite{BB1}, \cite{BB2})\footnote{In many interesting cases, the categories of Cartier crystals and unit $F$-crystals are in fact equivalent, see \cite{Sta}*{Theorem 2.13} and \cite{BB1}*{Theorem 5.15}. Cartier crystals are an analogue of right $\mathcal D$-modules, while unit $F$-crystals are an analogue of left $\mathcal D$-modules.} and perhaps also for the $\Lambda$-modules defined by B\"ockle and Pink in \cite{BoPi}.
\end{Remark}
 
\underline{\bf{Digression 1: Modules and comodules}} (see e.g. \cite{DNR}, \cite{HMS}).

If $k$ is a field and $A$ is a $k$-algebra, $\text{Hom}_k(A,k)$ is naturally a coalgebra  that we denote $A^{\vee}$. If $C$ is a $k$-coalgebra, then $C^{\vee}=\text{Hom}_k(C,k)$ is a $k$-algebra. 

If $A$ is {\it finitely dimensional} as a $k$-vector space, then there is an equivalence of categories \footnote{The equivalence holds also for $R$-algebras which are finitely generated and projective as an $R$-module, where $R$ is a commutative noetherian ring.}
\begin{eqnarray*}
  \left\{\text{f.g. left } A\text{-mods} \right\}  &\longrightarrow & \left\{\text{f.g. right } A^{\vee}\text{-comods}\right\}\ , \\
  M &\longmapsto & M^{\vee}
\end{eqnarray*}
but this fails if $A$ is infinite dimensional.

\begin{Example}(Arapura, \cite{Ara}*{Example 2.2.8}):
 Let $\mathcal A$ be the category of $\mathbb Z$-graded, finitely dimensional complex vector spaces, $F: \mathcal A \longrightarrow  \mathbb C-\mathfrak{fvs}$ the functor that forgets the grading. Denote $\text{End}_{\mathbb C}(F)$ the complex algebra of natural transformations $F \longrightarrow F$. Then $\mathcal A$ can be identified with the category of comodules over $\text{End}_{\mathbb C}(F)^{\vee}\cong\mathbb C[x,x^{-1}]$. But the category of modules over $\text{End}_{\mathbb C}(F)\cong\prod_{\mathbb Z}\mathbb C$ is bigger.
\end{Example}

However, one can prove that if $C$ is a $k$-coalgebra, then there is a full subcategory $\mathcal R$ of the category of left $C^{\vee}$-modules such that the functor
\begin{eqnarray*}
 \left\{\text{f.g. right } C\text{-comods}\right\} &\longrightarrow& \mathcal R \\
 M &\longmapsto& M^{\vee}
\end{eqnarray*}
is an equivalence of categories (see e.g. \cite{DNR}*{2.2} for proofs and details, see also \cite{AB}*{2.5}). In particular, this can be applied to the case $C=A(D,T)$ (where $D$ is a diagram and $T$ a representation). In this case, $C^{\vee}=\text{End}_{\mathbb C}(T)$, see \cite{Ara}*{Lemma 2.1.1}\footnote{But please be aware of the different notation used in \cite{Ara}.}.
\medskip

\underline{\bf{Digression 2: $\mathcal D$-modules on singular varieties}}

The usual definition of $\mathcal D$-modules on singular varieties is based on Kashiwara's equivalence\footnote{If $i: X \hookrightarrow Y$ is a closed embedding of smooth manifolds,  $i_{\ast}$ identifies the category of (quasi-coherent, coherent, holonomic) $\mathcal D_X$-modules with the category of (quasi-coherent, coherent, holonomic) $\mathcal D_Y$-modules with support on $X$.}. However, if one wants also {\it functoriality}, the issue is more delicate. There are at least two possibilities, both of them ``crystalline'':

\begin{enumerate}

\item One of them appears in \cite{QHE}*{7.10} and \cite{Gin}*{3.3.15} \footnote{The manuscript of Beilinson and Drinfeld is available at \\ https://www.math.uchicago.edu/$\sim$mitya/langlands/hitchin/BD-hitchin.pdf \\ Ginzburg's lecture notes are available at \\ www.math.harvard.edu/$\sim$gaitsgde/grad\_2009/Ginzburg.pdf}: If $f:X \longrightarrow X'$ is a finite morphism of schemes, denote $f^{(!)}: \mathcal O_{X'}-\mathfrak{mod} \longrightarrow \mathcal O_X-\mathfrak{mod}$ Grothendieck's functor (see \cite{QHE}*{7.10.1}). A closed embedding of schemes $Y \hookrightarrow X$ is a nilpotent thickening of $Y$ if $Y$ is defined by a nilpotent ideal. 

Then, a $\mathcal D$-crystal on $Y$ is a family $\{F_X\}_X$ of $\mathcal O_X$-modules, where $X$ runs over the nilpotent thickenings of $Y$, such that for any $f:X \longrightarrow X'$ making the diagram
  \[
    \begin{tikzcd}
    {} & Y \arrow{dr} \arrow{dl}\\
    X  \arrow{rr}{f} && X'
    \end{tikzcd}
   \]
commutative, there is an isomorphism $f^{(!)}(F_{X'})\cong F_X$, and these isomorphisms are compatible with composition. $\mathcal D$-crystals form a category, and direct and inverse images can be defined. If $Y$ is smooth, there is an equivalence between the category of right $\mathcal D_Y$-modules and the category of $\mathcal D$-crystals on $Y$. In fact this holds also in the non-smooth case, if $\mathcal D_Y$-modules are defined by taking an embedding of $Y$ in a smooth scheme, see \cite{Gin} and \cite{QHE} for details. 

\item Another approach is due to M. Saito (see \cite{Sai}).
Saito works in the analytic category, given an analytic space $Y$, he considers all closed embeddings $i: U \hookrightarrow Z$ where $U$ is open in $Y$ and $Z$ is a smooth variety. A $\mathcal D$-module on $Y$ is a compatible choice of a $\mathcal D_Z$-module on each such $Z$. 
He gives details on functoriality, proves a Riemann-Hilbert correspondence, etc...
\end{enumerate}

\bibliographystyle{amsplain}
\bibliography{expoenero.bib}
\end{document}